\documentclass[11pt]{article}
\usepackage{amssymb}
\newtheorem{Theorem}{Theorem}[section]
\newtheorem{Theorem--Definition}[Theorem]{Theorem--Definition}

\newtheorem{Lemma--Definition}[Theorem]{Lemma--Definition}
\newtheorem{Proposition}[Theorem]{Proposition}
\newtheorem{Proposition--Definition}[Theorem]{Proposition--Definition}

\newtheorem{Remark--Definition}[Theorem]{Remark--Definition}%
\newtheorem{Definition}[Theorem]{Definition}

\title{Coarsifications of the module isomorphism relation}
\author{Peteris Daugulis\\Department of Mathematics, Daugavpils University\\ Parades 1, Daugavpils, Latvia }
\date{February 1, 2003}
\begin{document}
\maketitle \pagenumbering{arabic}
\begin{abstract}
In this paper new equivalence relations on the category $Mod(A)$
for any associative algebra $A$ and several related results are
given. The new equivalence relations are defined using
restrictions to subalgebras and the action of algebra
automorphisms on modules.
\medskip

MSC 16G30, 16G60.
\end{abstract}
\section{Introduction}
\indent The classic representation/module theory of associative
algebras studies the module categories over algebras with their
algebra structures having been fixed once and for all. Modules are
usually studied up to isomorphism and isomorphism classes are
orbits of the set of objects of $Mod(A)$ under the conjugation
action of the general linear group of the underlying vector space
of the module. In this paper we introduce new equivalence
relations in the module category which are strictly coarser than
the isomorphism relation. This is motivated by a desire to make
progress in the classification of indecomposable modules in the
difficult cases such as the algebras of wild representation type.

There are two main ideas. The first idea is to consider
restrictions of a module to all subalgebras and collect the
information into a single structure. This idea is not entirely new
since restrictions to subalgebras of finite representation type
have been studied for many years. The second idea is   based on an
additional action on the category of modules - the action of the
automorphism group of the algebra. The orbits of this action
define an equivalence relation which is coarser than the standard
isomorphism equivalence relation. This equivalence relation allows
us to study representations of the algebra up to algebra
automorphisms. The described new equivalences relations can
contribute to the understanding of classic module categories in
the tame and wild cases, in particular, some parameters of
multiparameter families of indecomposable modules can be
attributed to the algebra automorphism action etc.

In this introductory paper we show give several easy first results
and examples for finite dimensional algebras.

\section{Definitions}

\subsubsection{Restriction isomomorphism}

\indent Let $k$ be a field, $A$ an associative  $k$-algebra. All
modules considered in this paper are left modules. We denote the
isomorphism class of a module $M$ by $[M]$. If the algebra $A$ is
given by a presentation $\langle X|\rho \rangle $ then the image
of $X$ under the algebra homomorphism defining the $A$-module $M$
is called the set of $A$-action generators for $M$. The categories
of all and finitely generated $A$-modules are denoted by $Mod(A)$
and $mod(A)$, respectively. By a subalgebra we mean a subalgebra
with the multiplicative identity. Given a subalgebra $A'\subseteq
A$ we denote the restriction of $M$ to $A'$ by $res_{A'}(M)$. By
$e_{ij}$ we denote the basis matrix having one noonzero entry at
position $(i,j)$.

\begin{Definition} $A$-modules $M_{1}$ and $M_{2}$ are called
restriction isomorphic (R-isomorphic, denoted by $[M_{1}]_{R}=
[M_{2}]_{R}$) provided for every proper subalgebra $A'\subset A$

\begin{equation} [res_{A'}(M_{1})] = [res_{A'}(M_{2})]
\end{equation}
\end{Definition}

\begin{Definition} Given a $A$-module $M$ and a proper subalgebra $A'\subseteq A$ the function $$A'\mapsto [res_{A'}(M)]$$
 is called
the restriction function (total restriction function) of $M$.

\end{Definition}

\begin{Definition}  $A$-module $M$ is called $R$-decomposable
provided $res_{A'}(M)$ is decomposable for any maximal proper
subalgebra $A'\subset A$.
\end{Definition}

\begin{Definition} $A$-modules $M_{1}$ and $M_{2}$ are called
restriction distinct ($R$-distinct) provided for every proper
subalgebra $A'\subset A$
\begin{equation} [res_{A'}(M_{1})] \neq [res_{A'}(M_{2})]
\end{equation}
\end{Definition}

\subsection{Action of algebra automorphism and twisted isomorphism}

\begin{Definition} $A$-modules $M_{1}$ and $M_{2}$ are called
twisted isomorphic (T-isomorphic, denoted by $[M_{1}]_{T}=
[M_{2}]_{T}$) provided there exists an $A$-automorphism
$f:A\rightarrow A$ and a $k$-linear isomorphism (twisted
isomorphism or $T$-isomorphism) $\varphi :M_{1} \rightarrow M_{2}$
such that for every $a\in A$

\begin{equation} f(a)\circ \varphi =\varphi  \circ a
\end{equation}
\end{Definition}

\begin{Definition} Let $f\in Aut(A)$, $M$ -
a $A$-module with the list of $A$-action generators
$(x_{1},...,x_{n})$. $f(M)$ (M twisted by f) is defined as the
module with the same underlying space and the action generators
$(f(x_{1}),...,f(x_{n}))$. We call the set $\bigcup\limits_{f\in
Aut(A)} f(M)$ the $T$-orbit of $M$.
\end{Definition}

\noindent We note that T-orbits of $A$-module isomorphism classes
are just equivalence classes under the T-isomorphism equivalence
relation and $[M_{1}]_{T}= [M_{2}]_{T}$ iff $[M_{1}]= [f(M_{2})]$
for some $f\in Aut(A)$.

\begin{Definition}  We call two  $A$-modules $M_{1}$  and  $M_{2}$ restriction-twisted
isomorphic (RT-isomorphic) provided  all restrictions to maximal
proper subalgebras are T-isomorphic.

\end{Definition}

\section{Results related to the $R$-isomorphism}

\subsection{$R$-isomorphism is strictly coarser}

\begin{Proposition} $R$-isomorphism is an equivalence relation on
$Obj(Mod(A))$ which is coars\-er than the standard module
isomorphism equivalence relation. Both restriction functions are
module invariants.
\end{Proposition}

Sets of module restriction isomorphism types are already being
studied in special cases such as modules over modular group
algebras restricted to subalgebras isomorphic to modular group
algebras for elementary abelian groups of rank $1$.

\begin{Proposition} For algebras of tame and wild representation type
$R$-isomorphism is strictly coarser than the standard module
isomorphism equivalence relation.
\end{Proposition}
\noindent \textbf{Proof} We exhibit two counterexamples. For the
tame case take $A=k[X,Y]/(X,Y)^{2}$. Let $M_{1}$ be the
indecomposable $A$-module of dimension 3 with
$dim_{k}(Soc(M_{1}))=1$ and $M_{2}$ be the indecomposable
$A$-module of dimension 3 with $dim_{k}(Soc(M_{2}))=2$. All
restrictions for both $M_{1}$ and $M_{2}$ to maximal proper
subalgebras decompose as direct sums of the trivial module and the
indecomposable two dimensional module therefore $[M_{1}]_{R}=
[M_{2}]_{R}$. For the wild case let $A=k[X,Y,Z]/(X,Y,Z)^{2}$. Let
$M_{1}$ and $M_{2}$ be modules of dimension 6 given by $(X,Y,Z)$
action matrix triples $(x,y,z)$ and $(x,z,y)$, respectively, where

\begin{equation}
\left\{%
\begin{array}{ll}
    x=e_{41}+e_{52}+e_{63}, \\
    y=e_{42}, \\
    z=e_{53}. \\
\end{array}%
\right.
\end{equation}

\noindent By direct matrix computations it can be proved that
$[M_{1}] \neq [M_{2}]$ and $[M_{1}]_{R}=[ M_{2}]_{R}$.

\subsection{Indecomposable $R$-decomposable modules}

\begin{Proposition} There exist algebras of  wild representation
type with indecomposable $R$-decomposable modules.
\end{Proposition}
\noindent \textbf{Proof} We exhibit an example. Let
$A=k[X,Y,Z]/(X,Y,Z)^{2}$. Let $M$ be the module of dimension $4$
given by $(X,Y,Z)$ action matrix triple $(x,y,z)$  where

\begin{equation}
\left\{%
\begin{array}{ll}
    x=e_{31}, \\
    y=e_{32}, \\
   z=e_{31}+e_{42}. \\
\end{array}%
\right.
\end{equation}

\noindent By direct computation it can be shown that $M$ is
indecomposable and $R$-decomposable.

\subsection{$R$-distinct modules}

\begin{Proposition} There exist algebras of  wild representation
type with $R$-distinct modules.
\end{Proposition}
\noindent \textbf{Proof} Let $A=k[X,Y,Z]/(X,Y,Z)^{2}$. Let $M_{1}$
and $M_{2}$ be modules of dimension $4$ given by $(X,Y,Z)$ action
matrix triples $(x_{i},y_{i},z_{i})$ where

\begin{equation}\left\{%
\begin{array}{ll}
    x_{1}=0, \\y_{1}=e_{41},\\ z_{1}=e_{31}+e_{42}\\
\end{array}%
\right.
\end{equation}
and
\begin{equation}\left\{%
\begin{array}{ll}
    x_{2}=e_{41}, \\y_{2}=e_{31}+e_{42},\\ z_{2}=e_{42}.\\
\end{array}%
\right.
\end{equation}

\noindent By direct matrix computations it can be shown that
$M_{1}$ and $M_{2}$ are $R$-distinct.

\section{Results related to the $T$-isomorphism}

\subsection{$T$-isomorphism ir strictly coarser}

\begin{Proposition} $T$-isomorphism is an equivalence relation on
$Obj(Mod(A))$ which is coars\-er than the standard module
isomorphism equivalence relation.
\end{Proposition}

\noindent \textbf{Proof} Equivalence property of the
$T$-isomorphism relation follows from the group properties of
$Aut(A)$. The strictness is proved below in the next subsection.

\subsection{Group action property and preservation of indecomposibility}

\begin{Proposition} The map $Obj(mod(A)) \rightarrow Obj(mod(A))$ defined by $M \mapsto f(M)$
is a group action which preserves indecomposability.
\end{Proposition}
\noindent \textbf{Proof} The finitely generated $A$-modules $M$
and $f(M)$ have the same endomorphism algebras therefore they are
simultaneously indecomposable or decomposable.

\subsection{Comparing $T$-isomorphism and $R$-isomorphism}

\begin{Proposition} $T$-isomorphism does not imply $R$-isomorphism
and vice versa.

\end{Proposition}
\noindent \textbf{Proof} Counterexamples for the first part of the
statement are easy to find for the algebra $A=k[X,Y]/(X,Y)^{2}$.

To exhibit a counterexample for the second part ($R$-isomorphism
does not imply $T$-isomorphism) take $A=k[X,Y]/(X,Y,Z)^{2}$. Let
$M_{1}$ and $M_{2}$ be modules of dimension $6$ given by $(X,Y,Z)$
action matrix triples $(x_{i},y_{i},z_{i})$ where

\begin{equation}
\left\{%
\begin{array}{ll}
   x_{1}=e_{51}+e_{42},\\ y_{1}=e_{61}+e_{53},\\
z_{1}=e_{52}+e_{43}+e_{63}.\\
\end{array}%
\right.
\end{equation}
and
\begin{equation}
\left\{%
\begin{array}{ll}
     x_{2}=x_{1}, \\y_{2}=y_{1},\\ z_{2}=e_{41}+e_{62}+e_{63}.\\
\end{array}%
\right.
\end{equation}

\noindent By direct matrix computations it can be shown that
$M_{1}$ and $M_{2}$ are $R$-isomorphic but not $T$-isomorphic.

\subsection{One-parameter families of indecomposable modules over
tame algebras as T-orbits}

We assume that $k$ is algebraically closed and all $A$-modules
considered below are finitely generated. By a family of
indecomposable modules we mean a subset of isomorphism classes of
$Obj(mod(A))$ continuously depending on arguments taking values in
an open subset of $k^{r}$ for some $r\in \mathbb{N}$, $r$ is said
to be the number of parameters of the family. In this section we
show that for several classic tame algebras one-parameter families
of indecomposable finitely generated modules are $T$-orbits or
their subsets. This also shows that $T$-isomorphism is strictly
coarser than the standard isomorphism.

\subsubsection{$A=k[X]$} Modules in one-parameter families of
indecomposable finitely generated $A$-modules are isomorphic to
modules with $X$ action matrices being Jordan block matrices.
%
\noindent The group $Aut(A)$ is the group of $A$-endomaps which
fix the element $1\in A$  and which are obtained by extending the
maps of form $X\mapsto aX+b$, $a\ne 0$ according to algebra
structure of $A$. For any $\lambda_{1}\in k$,$\lambda_{2}\in k$,
$\lambda_{1}\ne \lambda_{2}$ and any $a'\ne 0$ we can find $b'\in
k$ such that $\lambda_{2}=a'\lambda_{1}+b'$. Choosing the
$A$-automorphism $f$ generated by the map $X\mapsto a'X+b'$ and an
appropriate diagonal basis change $\varphi$ which makes
off-diagonal elements equal to 1 we see that for any $n\in
\mathbb{N}$ we have $[J(\lambda_{1},n)]_{T}=
[J(\lambda_{2},n)]_{T}$. Thus we have the following proposition.

\begin{Proposition} For any $n\in \mathbb{N}$ the one-parameter family
$\{J(\lambda,n)\}_{\lambda \in k}$ is a $T$-orbit.
\end{Proposition}

\subsubsection{$A=k[X,Y]/(X,Y)^{2}$} Modules in one-parameter
families of indecomposable finitely generated $A$-modules are
isomorphic to modules in the families
$\{K(\lambda,n)\}_{\lambda\in k\bigcup \lambda =\infty}$ where
modules $K(\lambda,n)$ are given by action generators
$K_{X}(\lambda,n)$ and $K_{Y}(\lambda,n)$ for $X$ and $Y$,
respectively, which can be chosen as

\begin{equation}
\left\{%
\begin{array}{ll}
   K_{X}(\lambda,n)=\left[ \begin{array}{c|c}
O_{n}&O_{n}\\
\hline E_{n}&O_{n}
\end{array} \right ]\\
,K_{Y}(\lambda,n)=\left[ \begin{array}{c|c}
O_{n}&O_{n}\\
\hline J(\lambda,n)& O_n
\end{array} \right ]\\
\end{array}%
\right.
\end{equation}

\noindent  for any $\lambda \in k$ and additionally

\begin{equation}\left\{%
\begin{array}{ll}
    K_{X}(\infty,n)=\left[ \begin{array}{c|c}
O_{n}&O_{n}\\
\hline J(0,n)&O_{n}
\end{array} \right ],\\
K_{Y}(\infty,n)=\left[ \begin{array}{c|c}
O_{n}&O_{n}\\
\hline E_{n}& O_n
\end{array} \right ]\\
\end{array}%
\right.
\end{equation}

\noindent For more details about finitely generated $A$-modules
see \cite{bashev}. $Aut(A)\simeq GL(2,k)$ and is formed by the
maps which fix $1\in A$ and map $(X,Y)$ to
$(a_{11}X+a_{12}Y,a_{21}X+a_{22}Y)$ with $\det([a_{ij}])\ne 0$.
For any $\lambda_{1}\in k$, $\lambda_{2}\in k$, $\lambda_{i}\ne
\infty$, $\lambda_{1}\ne \lambda_{2}$ as in the previous example
we can find $b'\in k$ such that $\lambda_{2}=a'\lambda_{1}+b'$,
$a'\ne 0$. By taking the automorphism given by $(X,Y)\mapsto
(X,a'Y+b')$ and a diagonal basis change as in the previous example
we see that $[K(\lambda_{1},n)]_{T}= [K(\lambda_{2},n)]_{T}$. By
taking $\varphi =id$ and the $A$-automorphism permuting $X$ and
$Y$ we see that $[K(0,n)]_{T}= [K(\infty,n)]_{T}$.

\begin{Proposition}
For any $n\in \mathbb{N}$ and any $\lambda \in k$ the one
parameter family $\{K(\lambda,n\}_{\lambda \in K \bigcup \infty}$
is a $T$-orbit.
\end{Proposition}

We end considering this example by noting that $A$-module
isomorphism classes which do not lie in infinite families (the
indecomposable $A$-modules of odd dimension) coincide with
$T$-isomorphism classes.

\subsubsection{$A=k[X,Y]/(X^2,Y^2,(XY)^kX^{\epsilon_{1}},(YX)^kY^{\epsilon_{2}})$}

In this case $A$ is closely related to group algebras for dihedral
2-groups over fields of characteristic 2. Modules in one-parameter
families of indecomposable finitely generated $A$-modules are
isomorphic to module in the families $\{B(w,\lambda,m)\}_{\lambda
\in k,\lambda\ne 0}$ where modules $B(w,\lambda,m)$ are given by
action generators $B_{X}(w,\lambda,m)$ and $B_{Y}(w,\lambda,m)$
(of size $nm$) for $X$ and $Y$, respectively, which can be chosen
as

$$
B_{X}(w,\lambda,m)=\underbrace{B_{X}(w,\lambda)\oplus ... \oplus
B_{X}(w,\lambda)}_{m\ times}
$$

\begin{equation}
B_{Y}(w,\lambda,m)=\left[ \begin{array}{c|c|c|c|c}
B_{Y}(w,\lambda)&O_{n}&O_{n}&...&O_{n}\\
\hline
O_{n}+e_{nn}&B_{Y}(w,\lambda)&O_{n}&...&O_{n}\\
\hline
O_{n}&O_{n}+e_{nn}&B_{Y}(w,\lambda)&...&O_{n}\\
\hline
...&...&...&...&...\\
\hline
O_{n}&O_{n}&O_{n}&...&B_{Y}(w,\lambda)
\end{array} \right ]
\end{equation}

\noindent and the matrices $B_{X}(w,\lambda)$,$B_{Y}(w,\lambda)$
of size $n$ are determined by the oriented edge labelled graph
isomorphism class of an admissible oriented edge labelled cycle
$w$ and parameter $\lambda \in k$ (see \cite{ringel1}). An
admissible oriented edge labelled cycle has edges labelled by $X$
and $Y$, it is subject to conditions determined by the structure
of $A$ (see \cite{ringel1}).  Note that there is a partition of
the vector space basis into groups of $n$ elements corresponding
to subquotients $B_{1},...,B_{m}$ each of which is isomorphic to
the $A$-module $B(w,\lambda)=B(w,\lambda,1)$ given by action
matrices $B_{X}(w,\lambda)$ and $B_{Y}(w,\lambda)$. Note that
$B_{X}(w,\lambda)$ does not depend on $\lambda$ and
$B_{Y}(w,\lambda)$ has one of its nonzero entries equal to
$\lambda$ and other nonzero entries equal to $1$.

There are $A$-automorphisms which fix $1\in A$ and map $(X,Y)$ to
$(X,aY)$, we denote such an $A$-automorphism by $f_{a}$. We want
to show that for any one-parameter family
$\{B(w,\lambda,m)\}_{\lambda \in k}$ any two distinct elements
$B(w,\lambda_{1},m)$ and $B(w,\lambda_{2},m)$, $\lambda_{i}\ne
0$,$\lambda_{1}\ne \lambda_{2}$, are $T$-isomorphic via a suitable
$A$-automorphism of type $f_{a}$. We first show that
$[B(w,\lambda_{1})]_{T}=[B(w,\lambda_{2})]_{T}$. We start with
module $B(w,\lambda_{1})$ and apply an automorphism $f_{a}$ for
some $a$ to it. After this operation the action generator for $X$
does not change, the action generator of $Y$ is multiplied by $a$.
By consecutive diagonal basis changes we can isomorphically
transform $f_{a}(B(w,\lambda_{1}))$ to a module
$B(w,\lambda_{1}a^z)$ for some $z\in \mathbb{Z}$. Since $k$ is
algebraically closed we can solve the equation
$\lambda_{1}a^z=\lambda_{2}$ with respect to $a$ and thus we have
that $[B(w,\lambda_{1})]_{T}=[B(w,\lambda_{2})]_{T}$ via the
computed $A$-automorphism $f_{a}$. Now suppose we have two modules
$B(w,\lambda_{1},m)$ and $B(w,\lambda_{2},m)$ with $m>1$. We
produce a suitable twisted isomorphism $\varphi$ in two steps of
basis changes. We first transform separately each of the
subquotients $B_{i}$ as described above using diagonal basis
changes to have each subquotient isomorphic to $B(w,\lambda_{2})$.
In the second step we consecutively multiply all basis elements of
each subquotient $B_{i}$ by the same coefficient (which will be
different for different subquotients) to make the elements in the
off-diagonal blocks for the $Y$ matrix equal to 1. Thus we have
the following proposition.

\begin{Proposition}
For any $m\in \mathbb{N}$ and any admissible word $w$ the
one-parameter family $\{B(w,\lambda,m\}_{\lambda \in k,\lambda\ne
0}$ is a subset of a $T$-orbit.
\end{Proposition}

\subsection{One parameter family does not necessarily belong to one $T$-orbit}

Let $A=k[x,y]/(x^2,y^3,y^2-(xy)x).$ If $char(k)=2$ then $A$ is a
semidihedral group algebra.

\begin{Proposition} There are one parameter families of
indecomposable $A$-modules which do not belong to one $T$-orbit.
\end{Proposition}

\noindent \textbf{Proof}  We will produce an example. Let
$f:A\rightarrow A$ be an $A$-automorphism, $
    f(x)=\tilde{x},$
$f(y)=\tilde{y}$. We must have
$$
\left\{%
\begin{array}{ll}
    \tilde{x}^2=0,\\
    \tilde{y}^3=0,\\   \tilde{y}^2=\tilde{x}\tilde{y}\tilde{x}.\\
\end{array}%
\right.
$$

Let
$$\tilde{x}=a_{1}x+a_{2}y+a_{3}xy+a_{4}yx+a_{5}xyx+a_{6}yxy.$$

Relation $\tilde{x}^2=0$ implies $a_{2}=0$. For $f$ to be
bijective  we must have $a_{1}\ne 0$.
\bigskip

Let $M_{1}$ and $M_{2}$ be modules of dimension $2$ given by
$(X,Y)$ action matrix pairs $(x_{i},y_{i})$ where

\begin{equation}
\left\{%
\begin{array}{ll}
    x_{1}=e_{21},\\ y_{1}=0\\
\end{array}%
\right.
\end{equation}
and
\begin{equation} \left\{%
\begin{array}{ll}
     x_{2}=0, \\y_{2}=e_{21}.\\
\end{array}%
\right.
\end{equation}

$M_{1}$ and $M_{2}$ belong to a one parameter family of
indecomposable $A$-modules defined by matrix pairs
$\{(e_{21},\lambda e_{21})\}\cup \{(0,e_{21})\}$.

If $M_{1}$ and $M_{2}$ would belong to one $T$-orbit then then
there would exist an $A$-automorphism $f$ such that
$f(x_{1})=x_{2}=0$. This is not possible since $f$ coefficient at
$x$ is not $0$.

\subsection{Examples of multiparameter module families over a
wild algebra as a subset of a $T$-orbit}

\subsubsection{A $2$-parameter family}

Let $A=k[X,Y,Z]/(X,Y,Z)^2$. $A$ is known to be a minimal wild
algebra, i.e. all its subalgebras are of tame or finite
representation type. Clearly $Aut(A)\simeq GL(3,k)$ and is formed
by maps which fix $1\in A$ and map $(X,Y,Z)$ to
$(a_{11}X+a_{12}Y+a_{13}Z$,$a_{21}X+a_{22}Y+a_{23}Z$,$a_{31}X+a_{32}Y+a_{33}Z)$
with $det(a_{ij})\ne 0$. Consider a family
$\mathcal{C}_{2}=\{C(\alpha,\beta)\}$, $\alpha\in k$,$\beta\in
k$,$\alpha\ne 0$,$\beta\ne 0$, of nonisomorphic indecomposable
$A$-modules of dimension 2 given by action generators

\begin{equation}
\left\{%
\begin{array}{ll}
    C_{X}(\alpha,\beta)=\left[ \begin{array}{cc}
0&0\\
1&0
\end{array} \right ]
,\\C_{Y}(\alpha,\beta)=\left[ \begin{array}{cc}
0&0\\
\alpha &0
\end{array} \right ]
,\\C_{Z}(\alpha,\beta)=\left[ \begin{array}{cc}
0&0\\
\beta &0
\end{array} \right ]\\
\end{array}%
\right.
\end{equation}

\noindent By considering diagonal $A$-automorphisms of form
$(X,Y,Z)\rightarrow (X,aY,bZ)$ we see that $[C(\alpha,\beta)]_{T}=
[C(1,1)]_{T}$ and thus the family $\mathcal{C}_{2}$ is a subset of
a $T$-orbit.

\subsubsection{A $3$-parameter family}

Let $A=k[X,Y,Z]/(X,Y,Z)^2$. Consider a family
$\mathcal{C}_{3}=\{C(\alpha,\beta,\gamma)\}$, $\alpha\in
k$,$\beta\in k$,$\alpha\ne 0$,$\beta\ne 0$, $\gamma\ne 0$, of
nonisomorphic indecomposable $A$-modules of dimension $5$ given by
action generators

\begin{equation}
\left\{%
\begin{array}{ll}
    x(\alpha,\beta)=e_{41}+e_{32} ,\\
y(\alpha,\beta)=e_{51}+e_{42} ,\\ z(\alpha,\beta)= \alpha
e_{32}+\beta e_{41}+\gamma (e_{51}+e_{42}).\\
\end{array}%
\right.
\end{equation}

By direct computation it can be shown that $\mathcal{C}_{3}$ is a
subset of a $T$-orbit.

\section{Research directions}

In this subsection we list several research problems related to
the described coarsifications of the isomorphism.

\subsection{Multiparameter families of indecomposable modules}

\begin{enumerate}

\item Do there exist multiparameter families of indecomposable
nonisomorphic and $X$-isomorphic modules ($X\in \{R,T,RT\}$)?
Describe such families.

\item Do there exist multiparameter families of indecomposable
nonisomorphic modules which are transversal with respect to the
$X$-isomorphism orbits (all elements are pairwise
non-$X$-isomorphic)? Describe such families.

\end{enumerate}

\subsection{Restriction functions}

\begin{enumerate}

\item What are the admissible restriction functions?

\item What are the sets of nonisomorphic modules for each value of
the restriction function? Describe the distribution. In
particular, for which values of the restriction function the sets
of nonisomorphic modules are wild (tame, finite)?

\item Describe restriction functions of modules in almost split
sequences.

\item Describe $R$-decomposable modules.

\end{enumerate}

\subsection{$T$-orbits}

\begin{enumerate}

\item What are the admissible $T$-orbits?

\item How do $T$-orbits depends on module structure of its
representatives?

\item Describe multiparameter families of indecomposable modules
which belong (do not belong) to a single $T$-orbit.

\end{enumerate}


\begin{thebibliography}{100}
\addcontentsline{toc}{chapter}{Bibliography}
\bibitem{bashev} V.A.Bashev, Representations of the group $Z_{2}\times Z_{2}$ in a field of characteristic 2 (Russian), Dokl.Akad.Nauk SSSR 141 (1961),1015-1018.


\bibitem{ringel1} C.M.Ringel, The indecomposable representations of the dihedral 2-ggroups, Math.Ann.214 (1975),19-34.

\end{thebibliography}
\end{document}